\documentclass[12pt, 14paper,reqno]{amsart}
\usepackage{amsmath,amsfonts,amssymb}
\usepackage{mathrsfs}
\usepackage{longtable}
\usepackage[breaklinks]{hyperref}
\usepackage{graphicx}
\usepackage{xcolor}

\theoremstyle{plain}
\newtheorem{thm}{Theorem}[section]
\newtheorem{lem}{Lemma}[section]

\newtheorem{thma}{Theorem}

\theoremstyle{proof}
\makeatletter
\@namedef{subjclassname@2020}{%
  \textup{2020} Mathematics Subject Classification}
\makeatother

\numberwithin{equation}{section}

\begin{document} 
\title{On the exponential Diophantine equation $x^2+p^mq^n=2y^p$}
\author{Kalyan Chakraborty and Azizul Hoque}
\address{Kalyan Chakraborty @Kerala School of Mathematics, Kozhikode-673571, Kerala, India.}
\email{kalychak@ksom.res.in}
\address{Azizul Hoque @Department of Mathematics, Faculty of Science, Rangapara College, Rangapara, Sonitpur-784505, India.}
\email{ahoque.ms@gmail.com}
\date{\today}
\maketitle

\begin{abstract}
We study the exponential Diophantine equation $x^2+p^mq^n=2y^p$ in positive integers $x,y,m,n$, and odd primes $p$ and $q$ using primitive divisors of Lehmer sequences in combination with elementary number theory. We discuss the solvability of this equation.
\end{abstract}

\section{Introduction} 
The exponential Diophantine equation, so called generalized Lebesgue-Nagell equation 
\begin{equation}\label{eqs}
x^2+d=\lambda y^n,~~~ x, y, n \in \mathbb{N}, n\geq 3, \lambda=1, 2, 4,
\end{equation}
has a long and rich history. Its study, for specific values of $d, n$ and $\lambda$, dates back to the works of Fermat and Euler. It was Lebesgue who first proved that \eqref{eqs} has no solution when $(d, \lambda) = (1, 1)$.  Nowadays there are many results on the solutions of \eqref{eqs} for various values of $d$, and also for $d$ ranging in some infinite set of positive integers. A beautiful survey on this very interesting topic can be found in \cite{LS20}. A generalization of \eqref{eqs}, that is, the Diophantine equation,
\begin{equation}\label{eqs1} cx^2 + d^m = 4y^n, ~~~ x, y, n \in \mathbb{N}, n\geq 3,\end{equation}
for given square-free positive integers $c$ and $d$, has been studied by several authors in (cf. \cite{BI02, BU01, CHS20, CHS-21, DGS20, LE95}). 
In \cite{BI02}, Bilu studied \eqref{eqs1} for `odd $m$' under certain conditions. He pointed out that there is a flaw in \cite{LE95}  and a {\it fortiori} Bugeaud's result in \cite{BU01}, and  corrected that inaccuracy. In contrast, \eqref{eqs1} has been studied in \cite{CHS-21, DGS20} for `even $m$'. The authors completely solved \eqref{eqs1} in \cite{BHS-19, CHS20} when $c=1$ and $d=1,2,3,7,11,19,43,67,163$. Recently, the second author deeply investigated an analogous generalization of \eqref{eqs} in \cite{HO20} when $\lambda=2$. In \cite{PA21}, Patel deeply studied a similar generalization of \eqref{eqs} for $\lambda=1$. More precisely, she presented a practical method for solving the Diophantine equation, $$cx^2+d=y^n, ~~ x,y, n\in \mathbb{N},~ \gcd(cx,d,y)=1,~ n\geq 3,$$
satisfying certain conditions.        

Many authors considered \eqref{eqs} when $d$ is a product of primes with different exponents and investigate their solutions (cf. \cite{AZ20, CH22, CHS21, DA11, LSTZ15}). In particular,  Alan and Zengin \cite{AZ20} (res. Chakraborty et al. \cite{CHS21}) completely solved \eqref{eqs} when $\lambda=1$ and $d=3^a41^b$ (res. $d= 2^a11^b19^c$). In \cite{DA11, LSTZ15}, the equation \eqref{eqs} has been solved when $d=2^ap^b$ with odd prime $p$. Recently, Godinho and Neumann \cite{GN21} obtained some conditions for the existence of the solutions of the following variation: $$x^2+p^mq^n=y^p,~~~x,y, m, n\in \mathbb{N}, \gcd(x,y)=1,$$ where $p$ and $q$ are distinct odd primes with $p\geq 5$. They used primitive divisors of certain Lucas sequences as well as the factorization of certain polynomials.

 In this paper, we investigate the integer solutions $(x,y,m, n)$ of the slightly different variation of the above equation, viz. 
\begin{equation}\label{eqp1}
x^2+p^mq^n=2y^p,~~ x,y,m,n\in\mathbb{N},  \gcd(x,y)=1,
\end{equation}
where $p>3$ and $q\geq 3$ are distinct primes. This equation has been well investigated when either $m=0$ or $n=0$ (cf. \cite{TE07, ZLT12}). We consider both $m$ and $n$ to be positive integers and obtain the conditions for which \eqref{eqp1} has no solution. The main result is:
\begin{thm}\label{thm} Let $p>3$ and $q\geq 3$ be distinct primes. Assume that $d$ is the square-free part of $p^mq^n$ and $h(-d)$ denotes the class number of $\mathbb{Q}(\sqrt{-d})$. If $n$ is odd, $d\not\equiv 7\pmod 8$ and $\gcd(p, h(-d))=1$, then \eqref{eqp1} has no solution, except $(p,q,x,y,m,n)=(5,17,21417,47,3,1)$. Further for even $n$, \eqref{eqp1} has no solution when $q^{n/2}\not\equiv \pm 1\pmod p$ and one of the following holds:
\begin{itemize}
\item[(i)] $m$ is even,
\item[(ii)] $m$ is odd and $p\not\equiv 7\pmod 8 $.
\end{itemize}
\end{thm}
\subsection*{Remarks} The following comments are worth to be noted.
\begin{itemize}
\item[(i)] Let's assume that $N$ is a positive integer with $p\mid N$ and consider the following extension of \eqref{eqs}: 
$$x^2+p^mq^n=2y^N,~~ x,y,m,n, N\in\mathbb{N},  \gcd(x,y)=1.$$
Then Theorem \ref{thm} also holds for this equation too.   
\item[(ii)] Let $p\in\{13, 19,43,67,163\}$ and $q\in \{11,19,43,67,163\}$. For $p\ne q$,  \eqref{eqp1} has no solution. 
\item[(iii)] For the twin primes (both $p$ and $p+2$ are primes), then the Diophantine equation $$x^2+p^{2m}(p+2)^{2p}=2y^p,~~ x,y,m,n\in\mathbb{N},  \gcd(x,y)=1 $$ 
has no solution.
\end{itemize}

Our method largely relies on the prominent result of Bilu, Hanrot and Voutier \cite{BH01, VO95} concerning the primitive divisors of Lehmer sequences.

\section{Preliminary descent }
The following lemma is a part of  \cite[Corollary 3.1]{YU05}.
\begin{lem}\label{lemYU}
For a given square-free positive integer $d>3$, let $k\geq 3$ be an odd integer such that $\gcd(k, h(-d))=1$. Then all the positive integer solutions $(x,y,z)$ of the Diophantine equation 
\begin{equation*}\label{eqyu}
x^2+dy^2=2z^k,~~ \gcd(x,dy)=1,
\end{equation*} 
can be expressed as 
$$\frac{x+y\sqrt{-d}}{\sqrt{2}}=\varepsilon_1\left(\frac{a+\varepsilon_2 b\sqrt{-d}}{\sqrt{2}}\right)^k,$$
where $a$ and $b$ are positive integers satisfying $a^2+b^2d=2y$, $\gcd(a, bd)=1$, and $\varepsilon_1, \varepsilon_2\in\{-1, 1\}$.
\end{lem}
Before proceeding further we need to recall some definitions. A pair  of algebraic integers $\alpha$ and $\bar{\alpha}$ forms a Lehmer pair when $(\alpha +\bar{\alpha})^2$ and $\alpha\bar{\alpha}$ are two non-zero coprime rational integers, and $\alpha/\bar{\alpha}$ is not a root of unity. Given a positive integer $t$, the $t$-th Lehmer number corresponds to the pair $(\alpha, \bar{\alpha})$ is defined as 
$$\mathcal{L}_t(\alpha, \bar{\alpha})=\begin{cases}
\dfrac{\alpha^t-\bar{\alpha}^t}{\alpha-\bar{\alpha}} & \text{ when $t$ is odd,} \vspace{1mm}\\
\dfrac{\alpha^t-\bar{\alpha}^t}{\alpha^2-\bar{\alpha}^2} & \text{ when $t$ is even}.
\end{cases}$$
It is well known that all Lehmer numbers are non-zero rational integers.  The pairs $(\alpha_1, \bar{\alpha}_1)$ and $(\alpha_2, \bar{\alpha}_2)$ are called equivalent if $\alpha_1/\alpha_2=\bar{\alpha}_1/\bar{\alpha}_2\in \{\pm 1, \pm\sqrt{-1} \}$. A prime divisor $p$ of $\mathcal{L}_t(\alpha, \bar{\alpha})$ is primitive if $p\nmid(\alpha^2-\bar{\alpha}^2)^2
\mathcal{L}_1(\alpha, \bar{\alpha}) \mathcal{L}_2(\alpha, \bar{\alpha}) \cdots \mathcal{L}_{t-1}(\alpha, \bar{\alpha})$. 
 
 The following result about the primitive divisors of Lehmer numbers \cite[Theorem 1.4]{BH01} will be handy for us. 
\begin{thma}\label{thmBH}
For any integer $t>30$, the Lehmer numbers $\mathcal{L}_t (\alpha, \bar{\alpha}) $ have primitive divisors.
\end{thma}

Assume that  $(\alpha, \bar{\alpha})$ is a Lehmer pair.  Then we can write $\alpha=(\sqrt{a}\pm\sqrt{b})/2$ and $\bar{\alpha}=(\sqrt{a}\mp\sqrt{b})/2$ by taking $a=(\alpha+\bar{\alpha})^2$ and $b=(\alpha-\bar{\alpha})^2$. The pair $(a, b)$ is called the parameters corresponding to the Lehmer pair $(\alpha, \bar{\alpha})$. The next lemma  is extracted from \cite[Theorem 1]{VO95}.
\begin{lem}\label{lemVO}
Let $p$ be a prime such that $7\leq p\leq 29$. If the Lehmer numbers $\mathcal{L}_p(\alpha, \bar{\alpha})$ have no primitive divisor, then up to equivalence, the parameters $(a, b)$ of the corresponding Lehmer pair $(\alpha, \bar{\alpha})$ are given by: 
\begin{itemize}
\item[(i)] for $p=7$, $(a, b)=(1,-7), (1, -19), (3, -5), (5, -7), (13, -3), (14, -22)$;
\item[(ii)] for $p=13$, $ (a, b)=(1,-7)$.
\end{itemize}
\end{lem}
Let $F_k$ (resp. $L_k$) be the $k$-th term in the Fibonacci (resp. Lucas) sequence defined as follows:\\
$$ \begin{cases}
  F_{k+2}=F_k+F_{k+1}, ~~F_0=0,~~   F_1= 1,~~ k\geq 0,\\
 L_{k+2}=L_k+L_{k+1},~~ L_0=2,  ~~L_1=1, ~~k\geq 0.
\end{cases} 
$$
Then the following lemma follows from \cite[Theorem 1.3]{BH01}. 
\begin{lem}\label{lemBH}
Assume that the Lehmer numbers $\mathcal{L}_5(\alpha, \bar{\alpha})$ have no primitive divisor. Then up to equivalence, the parameters $(a, b)$ of the corresponding pair $(\alpha, \bar{\alpha})$ are given by 
$$(a, b)=\begin{cases}
(F_{k-2\varepsilon}, F_{k-2\varepsilon}-4F_k)\text{ with } k\geq 3,\\
 (L_{k-2\varepsilon}, L_{k-2\varepsilon}-4L_k)\text{ with } k\ne 1;
 \end{cases}
 $$
where $t\ne 0$ and $k\geq 0$ are any integers and $\varepsilon=\pm 1$.
\end{lem}
We also need some basic properties of Fibonacci numbers and Lucas numbers.  
\begin{thma}[{\cite[Theorems 2 and 4]{CO64}}]\label{thmCO} For an integer $k\geq 0$, let $F_k$ (resp. $L_k$) denote the $k$-th Fibonacci (resp. Lucas) number. Then
\begin{itemize}
\item[(i)] if $L_k=2x^2$, then $(k, x)=(0, 1), (6,3)$;
\item[(i)] if $F_k=2x^2$, then $(k,x)=(0,0), (3,1), (6,2)$.  
\end{itemize}
\end{thma}
It is not hard to prove the following lemma. 
\begin{lem}\label{lem51}
Assume that  $(k, d)=(7,5), ( 47, 85)$. Then the solutions of the Diophantine equation $$x^2+dz^2=2\times k^5,~~x, z\in \mathbb{N}$$
are given by 
$$(x,z )=\begin{cases} (63, 77), (147,49), (183,5)\hspace{2mm} \text{ if }(k, d)=(7,5);\\
(21417, 5)\hspace{30mm} \text{  if } (k, d)=(47,85).
\end{cases}
$$
\end{lem}

\section{Proof of Theorem \ref{thm}}
Assume that $(x,y,m,n)$ is a positive integer solution of \eqref{eqp1} for a given pair of distinct primes $p>3$ and $q\geq 3$. Then \eqref{eqp1} can be written as 
\begin{equation}\label{eqp2}
x^2+z^2d=2y^p, 
\end{equation}  
where,
\begin{align}\label{eqpx}
\begin{cases}
z&=p^{m_1}q^{n_1}\text{ with } m_1=\left \lfloor{\frac{m}{2}}\right \rfloor \text{ and }  n_1=\left \lfloor{\frac{n}{2}}\right \rfloor\\
d&=1,p,q,pq.
\end{cases}
\end{align}
Since both $p$ and $q$ are odd primes, so that $x$ is odd. If $y$ is even then  reading \eqref{eqp2} modulo $8$, we get $d\equiv 7\pmod 8$.  
This contradicts our assumption and thus $y$ is also odd. 

Note that by the assumption, $\gcd(p, h(-d))=1$ when $n$ is odd. Again for even $n$, we have $d=1$ when $m$ is even and $d=p$ when $m$ is odd. In the later case, it follows from the Dirichlet class number formula that $1\leq  h(-p) < p$ (see \cite[p. 67]{CO07-2}). Thus in either cases, $\gcd(p, h(-d))=1$.

Here, $\gcd(dz, x)=1$ as $\gcd(x, y)=1$. Since $\gcd(p, h(-d))=1$ with $p$ odd, so that by Lemma \ref{lemYU} we have 
\begin{equation}\label{eqp3}
\frac{x+z\sqrt{-d}}{\sqrt{2}}=\varepsilon_1\left(\frac{a+\varepsilon_2 b\sqrt{-d}}{\sqrt{2}}\right)^p,
\end{equation}
where $a$ and $b$ are positive coprime integers satisfying
 \begin{equation}\label{eqp4}
 a^2+b^2d=2y.
 \end{equation}
Here, $\varepsilon_1, \varepsilon_2$ are as defined in Lemma \ref{lemYU}.   Note that for $d=3$, $\varepsilon_1$ satisfies $\varepsilon_1^6=1$ and hence it can be absorbed into the $p$-th power as $p>3$.  


Assume that
\begin{equation}\label{eqp5}
\begin{cases}
\alpha=\dfrac{a+\varepsilon_2 b \sqrt{-d}}{\sqrt{2}},\\
\bar{\alpha}=\dfrac{a-\varepsilon_2 b\sqrt{-d}}{\sqrt{2}}.
\end{cases}
\end{equation}

By \eqref{eqp4}, we see that $\alpha$ satisfies the polynomial $X^4+2(y-a^2)X^2+y^2\in \mathbb{Z}[X]$, and hence $\alpha$ is an algebraic integer and so is $\bar{\alpha}$. Since $\gcd(a,bd)=1$, so that \eqref{eqp4} gives $\gcd(a, y)=1$, and thus $(\alpha+\bar{\alpha})^2=2a^2$ and $\alpha\bar{\alpha}=y$ are coprime (as $y$ is odd). 

Now we have, $$\frac{2a^2}{y}=\frac{(\alpha+\bar{\alpha})^2}{\alpha\bar{\alpha}}=\frac{\alpha}{\bar{\alpha}}+\frac{\bar{\alpha}}{\alpha}+2,$$ which implies
that $$y\left(\frac{\alpha}{\bar{\alpha}}\right)^2+2(y-a^2)\frac{\alpha}{\bar{\alpha}}+y=0.$$ 
Since $\gcd(2(y-a^2), y)=\gcd(2a,y)=1$, so that $\alpha/\bar{\alpha}$ is not an algebraic integer and thus it is not a root of unity. Therefore, $(\alpha, \bar{\alpha})$ is a Lehmer pair.  

As $p$ is odd prime, so that the corresponding Lehmer numbers are given by $$\mathcal{L}_p(\alpha, \bar{\alpha})=\frac{\alpha^p-\bar{\alpha}^p}{\alpha-\bar{\alpha}}.$$
We use \eqref{eqp2}, \eqref{eqp3} and \eqref{eqp5} to get 
\begin{equation}\label{eqp6}
|\mathcal{L}_p(\alpha, \bar{\alpha})|=\left| \frac{p^{m_1}q^{n_1}}{b}\right|.
\end{equation}
Since $\mathcal{L}_p(\alpha, \bar{\alpha})\in \mathbb{Z}$, so that $b\mid p^{m_1}q^{n_1}$ which ensures that  $p$ and  $q$ are the only candidates for the primitive divisors of $\mathcal{L}_p(\alpha, \bar{\alpha})$. However the fact `\textit{if $\ell$ is a primitive divisor of $\mathcal{L}_p(\alpha, \bar{\alpha})$, then $\ell\equiv \pm 1 \pmod p$}'  helps us to remove $p$ from the possibility of the primitive divisor. 

Now if $q$ is the primitive divisor of $\mathcal{L}_p(\alpha, \bar{\alpha})$ then $q\nmid (\alpha^2-\bar{\alpha}^2)^2$.  Here,
$(\alpha^2-\bar{\alpha}^2)^2=-4a^2b^2d $ and thus $q$ to be the primitive divisor of $\mathcal{L}_p(\alpha, \bar{\alpha})$ only when $n_1\geq 1$, $b=p^{m_2}$ and $d\in\{1,p\}$ with $0\leq m_2\leq m_1$. Therefore, we can conclude by \eqref{eqp6} that $\mathcal{L}_p(\alpha, \bar{\alpha})$ has no primitive divisor except for the cases $n_1\geq 1$, $d\in\{1,p\}$ and $b=p^{m_2}$ with $0\leq m_2\leq m_1$. Thus by Theorem \ref{thmBH}, there is no Lehmer number $\mathcal{L}_p(\alpha, \bar{\alpha})$ for $p>30$ and hence \eqref{eqp1} has no positive integer solution for $p>30$ except for $n_1\geq 1$, $d\in\{1,p\}$ and $b=p^{m_2}$ with $0\leq m_2\leq m_1$.
 
Since $\left((\alpha+\bar{\alpha})^2, (\alpha-\bar{\alpha})^2\right)=(2a^2,-2b^2d)$ is the parameters of the pair $(\alpha, \bar{\alpha})$, so that by Lemma \ref{lemVO} 
there is no Lehmer number $\mathcal{L}_p(\alpha, \bar{\alpha})$ for $p\geq 7$. Therefore \eqref{eqp1} has no positive integer solution for $p\geq 7$. 

Now for $p=5$,  Lemma \ref{lemBH} gives us 
\begin{itemize}
\item[(i)] $F_{k-2\varepsilon}=2a^2$,
\item[(ii)] $4F_k-F_{k-2\varepsilon}=2b^2d$ with $k\geq 3$,\\
or,
\item[(iii)] $L_{k-2\varepsilon}=2a^2$,
\item[(iv)] $4L_k-L_{k-2\varepsilon}=2b^2d$ with $k\ne 1$.
\end{itemize}
Applying Theorem \ref{thmCO} in (i), we get $(k, \varepsilon, a)\in \{(2,1,0),
 (5,1,1), (4,-1,2), (8,1,2)\}$. As $a$ is positive odd integer, so that $(k, \varepsilon, a)=(5,1,1)$ and thus by (ii), we have $4F_5-F_3=2b^2d$. This implies that $b^2d=9$, and hence by \eqref{eqp4}, we get $(y, d)=(5,1)$. Therefore \eqref{eqp2} becomes $x^2+z^2=2\times 5^5$, which gives $(x,z)=(3, 79), (79, 3)$ as $\gcd(x, y)=1$. These solutions lead to $m=0$, which is out of our consideration. Nevertheless, these solutions are listed in Table \ref{T1}. 
 
 As before,  by Theorem \ref{thmCO} we get $(k, \varepsilon, a)\in \{(4,-1,3), (8,1,3)\}$, and thus (ii) gives $(a,b, d)\in \{(3,1,5), (3,1,85)\}$. Therefore using \eqref{eqp4}, we get $(d, y)\in \{(5,7), (85, 47)\}$, and hence \eqref{eqp2} and Lemma \ref{lem51} together give us $(x,y,z, d)=(183, 7, 5, 5), (21417, 47, 5, 85)$. As $dz^2=p^mq^n$ (see \eqref{eqpx}), so that $(x,y,z, d)=(183, 7, 5, 5)$ implies that $n=0$ which is not possible as $n\geq 1$. Similarly, $(x,y,z, d)=(21417, 47, 5, 85)$ gives $(p, q, m, n)=(5, 17, 3, 1)$. Thus $(p, q, x, y, m, n)=(5, 17, 21417, 47, 3, 1)$, which  is a solution of \eqref{eqp1}. These solutions are listed in Table \ref{T1}.

{\small
\begin{longtable}{cccccc}
\caption{All the solutions of \eqref{eqp1}} \label{T1} \\

\hline \multicolumn{1}{c}{$x$} & \multicolumn{1}{c}{$y$} & \multicolumn{1}{c}{$p$}& \multicolumn{1}{c}{$q$} & \multicolumn{1}{c}{$m$} & \multicolumn{1}{c}{$n$} \\ \hline 
\endfirsthead
\endhead
\hline \multicolumn{6}{c}{{Continued on next page}} \\ \hline
\endfoot
\hline 
\endlastfoot
3&5&5&79&0&2\\ \hline
79&5&5&3&0&2\\ \hline
183&7&5&q&3&0\\ \hline
21417&47&5&17&3&1\\
\end{longtable}}

Now we consider the remaining case, that is $d\in\{1,p\}$ with $n_1\geq 1$ and $b=p^{m_2} ~~ (0\leq m_2\leq m_1)$. In the case, $|\mathcal{L}_p(\alpha, \bar{\alpha})|=p^{m_1-m_2}q^{n_1}$ with $n_1\geq 1$, which ensures that $q$ is the primitive divisor of $\mathcal{L}_p(\alpha, \bar{\alpha})$. Thus, we can not utilize the previous technique to find the positive integer solutions of \eqref{eqp1}. 

We equate the imaginary parts from both sides of \eqref{eqp3} to get the following:
$$2^{(p-1)/2}p^{m_1}q^{n_1}=\varepsilon_1\varepsilon_2b\sum\limits_{j=0}^{(p-1)/2}\binom{p}{2j+1}a^{p-2j-1}b^{2j}(-d)^j.
$$
Now $d\in\{1,p\}$ implies by \eqref{eqp2} that $n$ is even and thus $n_1=n/2$. Since $b=p^{m_2} ~~ (0\leq m_2\leq m_1)$, so that the above equation reduces to
\begin{equation}\label{eq10}
2^{(p-1)/2}p^{m_1-m_2}q^{n/2}=\varepsilon_1\varepsilon_2\sum\limits_{j=0}^{(p-1)/2}\binom{p}{2j+1}a^{p-2j-1}p^{2jm_2}(-d)^j.
\end{equation}

We first consider the case when $d=1$. In this case, $m$ is even and hence $m_1=m/2$. Thus \eqref{eq10} implies that 
\begin{equation}\label{eq11}
2^{(p-1)/2}p^{m/2-m_2}q^{n/2}=\varepsilon_1\varepsilon_2\sum\limits_{j=0}^{(p-1)/2}\binom{p}{2j+1}a^{p-2j-1}p^{2jm_2}(-1)^j.
\end{equation}
If $m_2=0$, then the above equation can be reduced to 
$$2^{(p-1)/2}p^{m/2}q^{n/2}\equiv \varepsilon_1\varepsilon_2(-1)^{(p-1)/2}\pmod p.$$
This implies that $m=0$, which contradicts to our assumption. Therefore $m_2\geq 1$ and hence reading \eqref{eq11} modulo $p^2$, we arrive at
$$2^{(p-1)/2} p^{m/2-m_2}q^{n/2}\equiv \varepsilon_1\varepsilon_2pa^{p-1}\pmod {p^2}.$$
This can be further reduced to 
$p^{m/2-m_2-1}q^{n/2}\equiv \pm 1\pmod p$, which implies that $m_2=m/2-1$ and $q^{n/2}\equiv \pm 1\pmod p$. This contradicts to the assumption.

Now we consider the final case, that is $d=p$. In this case, \eqref{eq10} implies that 
$$2^{(p-1)/2}p^{m_1-m_2}q^{n/2}=\varepsilon_1\varepsilon_2\sum\limits_{j=0}^{(p-1)/2}\binom{p}{2j+1}a^{p-2j-1}p^{2jm_2}(-p)^j.
$$
This further implies that 
$$2^{(p-1)/2}p^{m_1-m_2-1}q^{n/2}=\varepsilon_1\varepsilon_2\left( a^{p-1}-\sum\limits_{j=1}^{(p-1)/2}\binom{p}{2j+1}a^{p-2j-1}p^{2jm_2}(-p)^{j-1}\right).
$$
Since $\sum\limits_{j=1}^{(p-1)/2}\binom{p}{2j+1}a^{p-2j-1}p^{2jm_2}(-p)^{j-1}\equiv 0\pmod p$, so that the above equation implies that 
$$p^{m_1-m_2-1}q^{n/2}\equiv \pm 1\pmod p.$$
This further implies that $q^{n/2}\equiv \pm 1\pmod p$ with $m_1=m_2+1$, which contradicts our hypothesis. Thus we complete the proof.  

\subsection*{Acknowledgements}
The author is thankful to Professor  Andrzej Dabrowski  for providing a copy of \cite{DA11}. This work is supported by SERB MATRICS (Project No. MTR/2021/000762), Govt. of India. The authors are grateful to the anonymous referee for careful reading and valuable comments which have helped to improve this paper. 


\begin{thebibliography}{25}
\bibitem{AZ20} M. Alan and U. Zengin, {\it On the Diophantine equation $x^2 + 3^a41^b = y^n$}, Period. Math. Hung. {\bf 81} (2020), 284--291.

\bibitem{BHS-19} S. Bhatter, A. Hoque and R. Sharma, {\it On the solutions of a Lebesgue-Nagell type equation},  Acta Math. Hungar. {\bf 158} (2019), no. 1, 17--26.

\bibitem{BH01} Y. Bilu, G. Hanrot and P. M. Voutier, {\it Existence of primitive divisors of Lucas and Lehmer numbers (with an appendix by M. Mignotte)}, J. Reine Angew. Math. {\bf 539} (2001), 75--122.

\bibitem{BI02} Y. Bilu, {\it On Le's and Bugeaud's papers about the equation $ax^2 + b^{2m-1} = 4c^p$}, Monatsh. Math. {\bf 137} (2002), 1--3.

\bibitem{BU01} Y. Bugeaud, {\it On some exponential Diophantine equations},  Monatsh. Math. {\bf 132} (2001), 93--97.

\bibitem{CH22} K. Chakraborty and A. Hoque, {\it On the Diophantine equation $dx^2+p^{2a}q^{2b}=4y^p$}, Results Math. {\bf 77} (2022), no. 1, 11pp, article no. 18. 

\bibitem{CHS20} K. Chakraborty, A. Hoque and R. Sharma, {\it Complete solutions of certain Lebesgue-Ramanujan-Nagell equations}, Publ. Math. Debrecen {\bf 97} (2020), no. 3-4, 339--352.

\bibitem{CHS21} K. Chakraborty, A. Hoque and R. Sharma, {\it On the solutions of certain Lebesgue-Ramanujan-Nagell equations},  Rocky Mountain J. Math. {\bf 51} (2021), no. 2, 459--471.

\bibitem{CHS-21} K. Chakraborty, A. Hoque and K. Srinivas, {\it On the Diophantine equation $cx^2+p^{2m}=4y^n$}, Results Math. {\bf 76} (2021), no. 2, pp. 12, article no. 57.



\bibitem{CO07-2} H. Cohen, {\it Number Theory}, Volume II, Graduate Texts in Mathematics, Vol. 240, Springer, New York, 2007.

\bibitem{CO64}  J. H. E. Cohn, {\it Square Fibonacci numbers, etc.}, Fibonacci Quart. {\bf 2} 1964, no. 2, 109--113.

\bibitem{DA11} A. Dabrowski, {\it On the Lebesgue-Nagell equation}, Colloq. Math. {\bf 125} (2011), no. 2, 245--253.

\bibitem{DGS20} A. Dabrowski, N. G\"unhan and G. Soydan, {\it On a class of Lebesgue-Ljunggren-Nagell type equations}, J. Number Theory,  {\bf 215} (2020), 149--159.

\bibitem{HO20} A. Hoque, {\it On a class of Lebesgue-Ramanujan-Nagell equations}, Period. Math. Hungar. (to appear). 

\bibitem{GN21} H. Godinho and V. G. L. Neumann, {\it The Diophantine equation $x^2+p^aq^b=y^q$}, Int. J. Number Theory, {\bf 17} (2021), no. 09, 2113--2130.


\bibitem{LE95} M. Le,  {\it On the Diophantine equation $D_1x^2 + D_2 = 4y^n$},  Monatsh. Math. {\bf 120} (1995), 121--125.

\bibitem{LS20} M. Le and G. Soydan,  {\it A brief survey on the generalized Lebesgue-Ramanujan-Nagell  equation}, Surv. Math. Appl. {\bf 115} (2020), 473--523.

\bibitem{PA21}  V. Patel, {\it A Lucas-Lehmer approach to generalised Lebesgue-Ramanujan-Nagell equations}, Ramanujan J. {\bf 56} (2021), no.2,  
585--596.

\bibitem{TE07} Sz. Tengely, {\it On the Diophantine equation $x^2 + q^{2m} = 2y^p$}, Acta Arith. {\bf 127} (2007), no. 1, 71--86.

\bibitem{VO95} P. M. Voutier, {\it Primitive divisors of Lucas and Lehmer sequences}, Math. Comp. {\bf 64} (1995),
869--888.

\bibitem{YU05} P. Z. Yuan, {\it On the Diophantine equation $ax^2 + by^2 = ck^n$}, Indag. Math. (N. S.) {\bf 16} (2005), no. 2,
301--320.

\bibitem{ZLT12} H. Zhu, M. -H. Le and A. Togb\'{e}, {\it On the exponential Diophantine equation $x^2 + p^{2m} = 2y^n$}, Bull. Aust. Math. Soc. {\bf 86} (2012), 303--314.

\bibitem{LSTZ15} H. Zhu, M. Le, G. Soydan and A. Togb\'e, {\it On the exponential Diophantine equation $x^2 + 2^ap^b = y^n$}, Perod. Math. Hung. {\bf 70} (2015), 233--247.

\end{thebibliography}
\end{document}